\documentclass[a4paper,12pt]{article}
\usepackage[dvips]{epsfig}
\usepackage{amssymb}
\usepackage{amsmath}
\usepackage{epsfig}

\newcommand{\G}{\Gamma}

\newcommand{\la}{\langle}
\newcommand{\ra}{\rangle}
\newcommand{\qed}{\hfill\hbox{\rule{3pt}{6pt}} \medskip}
\newcommand{\proof}{{\sc Proof. }}
\newcommand{\wtg}{\widetilde{\G}_1}
\newcommand{\wtd}{\tilde{d}}

\newtheorem{theorem}{Theorem}[section]
\newtheorem{lemma}[theorem]{Lemma}
\newtheorem{corollary}[theorem]{Corollary}
\newtheorem{proposition}[theorem]{Proposition}

\newtheorem{convention}[theorem]{Convention}
\newtheorem{remark}[theorem]{Remark}

\newtheorem{question}[theorem]{Question}
\newtheorem{problem}[theorem]{Problem}
\newtheorem{construction}[theorem]{Construction}

\title{ON THE CONNECTIVITY OF BIPARTITE DISTANCE-BALANCED GRAPHS}

\author{\v{S}tefko Miklavi\v{c} \\ 
        UP PINT and UP FAMNIT \\
        University of Primorska \\
        Muzejski trg 2 \\ 
        6000 Koper, Slovenia \\
        stefko.miklavic@upr.si \and
        Primo\v{z} \v{S}parl \\
        Faculty of Education \\
        University of Ljubljana \\
        Kardeljeva plo\v{s}\v{c}ad 16 \\
        1000 Ljubljana, Slovenija \\
        primoz.sparl@pef.uni-lj.si}

\begin{document}
\maketitle

\begin{abstract}
A connected graph $\G$ is said to be {\it distance--balanced} whenever for any pair of adjacent vertices $u,v$ of
$\G$ the number of vertices closer to $u$ than to $v$ is equal to the number of
vertices closer to $v$ than to $u$. 
In [Bipartite graphs with balanced $(a,b)$-partitions, {\em Ars Combin.} {\bf 51} (1999), 113--119]
Handa asked whether every bipartite distance-balanced graph, that is not a cycle, is 3-connected. 
In this paper the Handa question is answered in the negative. Moreover, we show that a minimal
bipartite distance-balanced graph, that is not a cycle and is not 3-connected, has $18$ vertices and is unique.
In addition, we give a complete classification of non-$3$-connected bipartite distance-balanced graphs for which the minimal distance between
two vertices in a $2$-cut is three. All such graphs are regular and for each $k \geq 3$ there exists an infinite family of such graphs which are $k$-regular.

Furthermore, we determine a number of structural properties that a 
bipartite distance-balanced graph, which is not 3-connected, must have. As an application, 
we give a positive answer to the Handa question for the subfamily of
bipartite strongly distance-balanced graphs.
\end{abstract}

\section{Introduction}
\label{sec:intro}

Throughout this paper, all graphs are
connected, finite, undirected, without loops and multiple edges.
Given a graph $\G$ let $V(\G)$ and $E(\G)$ denote its vertex set and edge set, respectively.
For $u,v \in V(\G)$ we denote the distance between $u$ and $v$ by $d(u,v)$.
Furthermore, for any nonnegative integer $i$ and $u \in V(\G)$ let $N_i(u)=\{v \in V(\G) \mid d(u,v)=i\}$
(we abbreviate $N(u) = N_1(u)$).
For $W \subseteq V(\G)$ the subgraph of $\G$ induced by $W$ is denoted by $\la W \ra$
(we abbreviate $\G-W = \la V(\G) \setminus W \ra$). 
A {\em vertex cut} of $\G$ is a set $W \subseteq V(\G)$, such that $\G-W$ is disconnected. 
(A vertex cut of size $k$ is called a $k$-cut.) A graph is called $k$-{\em connected} if it has at least $k+1$ vertices and the size of the 
smallest vertex cut is at least $k$.

For any pair of vertices $u,v \in V(\G)$ we define $W_{uv}$ by
$$
  W_{uv} = \{ z \in V(\G) \mid d(u,z) < d(v,z) \}.
$$
A connected graph $\G$ is said to be {\it distance-balanced} whenever for any pair of adjacent vertices $u,v \in V(\G)$ we 
have
$$
 |W_{uv}| = |W_{vu}|.
$$

Distance-balanced graphs have been extensively studied, see \cite{BCPSSS, Ha, HN84, IKM, JKR, KMMM, KMMM2}.
It was shown in \cite[Lemma 2.1]{Ha} that every distance-balanced graph is $2$-connected. 
In the same paper (see also \cite{IKM}), Handa asked the following question.

\begin{question}[Handa \cite{Ha}]
\label{que}
Is every bipartite distance--balanced graph, that is not a cycle, 3-connected?
\end{question}

In \cite{Ha}, the positive answer to the above question was obtained for the family of
distance--balanced partial cubes. (In the last section of the paper we give an alternative proof of this result.) 
Motivated by the above question, we investigate the structural properties that a 
bipartite distance-balanced graph, which is not 3-connected and not isomorphic to a cycle, must have.
It turns out that the conditions on such graphs are quite restrictive. In fact, they enable us to answer the Handa question in the negative. 
Moreover, using these results we give a complete classification of non-$3$-connected bipartite distance-balanced graphs, for which the minimal distance between
the vertices in a $2$-cut is three. It turns out that all such graphs are $k$-regular for some $k \geq 3$, and that for every
such $k$ there is an infinite family of such graphs, one of order $2k\ell$ for each odd $\ell \geq 3$ (see Theorem~\ref{thm:main}).
We also show that the smallest non-$3$-connected bipartite distance-balanced graph that is not a cycle is unique and belongs to 
this family for $k = 3$ with $\ell = 3$ (its order is thus $18$). 
In contrast, we show that every bipartite strongly distance-balanced graph that is not a 
cycle is 3-connected (see Section \ref{sec:appl} for the formal definition of strongly distance-balanced graphs).

Throughout the paper we are using some of the results from \cite{Ha}. For the sake of self-containment we gather these 
results in the following proposition.

\begin{proposition}[\cite{Ha}]
\label{handa}
Let $\G$ be a distance-balanced graph with at least two edges. 
Then $\G$ is $2$-connected. Moreover, if $\G$ is bipartite, then the following hold.
\begin{itemize}
\item[(i)]
Let $x,y$ be vertices of $\G$ such that $\G - \{x,y\}$ is disconnected. Then $d(x,y) \ge 2$.
\item[(ii)]
If $\G$ is not a cycle, then the minimal degree of $\G$ is at least $3$.
\item[(iii)]
Assume $\G$ is not $3$-connected. Among all pairs of vertices $x,y \in V(\G)$ such that $\G - \{x,y\}$
is disconnected pick a pair $a,b$ for which $d(a,b)$ is minimal. 
Then $\G-\{a,b\}$ has exactly two components.
\item[(iv)]
Let $x,y$ be a pair of adjacent vertices of $\G$ and let  
$u \in W_{xy}$ and $v \in W_{yx}$. If $u$ and $v$ are adjacent, then $d(x,u)=d(y,v)$. 
\end{itemize}
\end{proposition}

For the rest of this paper we make the following assumptions and adopt the following notational convention.

\begin{figure}[h]
\begin{center}
\includegraphics[scale=0.6]{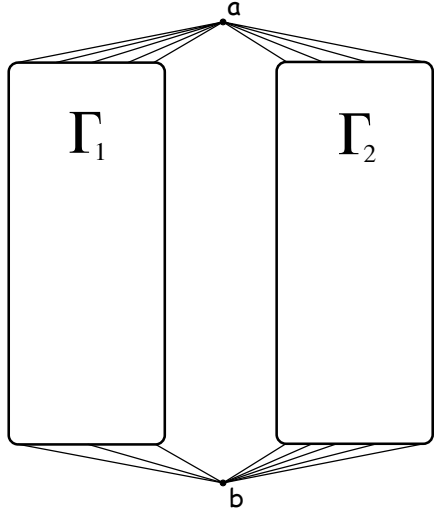}
\caption{The two components $\G_1$ and $\G_2$.}
\label{fig:gamma1and2}
\end{center}
\end{figure}

\begin{convention}
\label{blank}
{\em Let $\G$ denote a bipartite distance--balanced graph that is 
not a cycle and that is not 3-connected, and let $n$ denote its order. 
Since $\G$ is not $3$-connected it has a $2$-cut (recall that $\G$ is 2-connected by Proposition~ref{handa}).
Among all $2$-cuts pick a $2$-cut $\{a,b\}$ for which $d(a,b)$ is minimal. By Proposition~\ref{handa} we have that $d(a,b) \ge 2$ and $\G-\{a,b\}$ has exactly two components.
Denote these two components by $\G_1$ and $\G_2$ (see Figure~\ref{fig:gamma1and2}).  
Note that, since $\G$ is bipartite, we have $V(\G) = W_{uv} \cup W_{vu}$ for any pair of adjacent vertices $u,v \in V(\G)$.
In particular, $n$ is even and $|W_{uv}| = |W_{vu}| = n/2$.}
\end{convention}

\section{Good and bad vertices}
\label{sec:goodbad}

Let $\G$ be as in Convention \ref{blank}. In this section we introduce the notion of good and bad vertices of $\G$.
As we show below, one of the two components $\G_1$ and $\G_2$ coincides 
with the set of all bad vertices while the other, together with $a$ and $b$, coincides with the set of all good vertices.

\smallskip \noindent

A vertex $c \in V(\G)$ is called {\em good} if it lies on some shortest $a$-$b$ path of $\G$.
The vertices which are not good are called {\em bad}. First we make the following observation, 
which will be extensively used in the rest of the paper (without explicit reference to it). 

\begin{lemma}
\label{lem1}
Let $\G$ be as in Convention \ref{blank}. Pick $x \in \G_1$ and $y \in \G_2$.
Then $d(x,y) = \min\{d(x,a) + d(a,y), d(x,b) + d(b,y)\}$.
\end{lemma}
\proof
Since $\G-\{a,b\}$ is disconnected, every shortest path from $x$ to $y$ passes either through $a$ or through $b$. \qed

\begin{lemma}
\label{lem2}
Let $\G$ be as in Convention \ref{blank}.
If there is a bad vertex in $\G_i \; (i \in \{1,2\})$, then $|\G_i| \ge n/2$.
\end{lemma}
\proof
Without loss of generality we can assume that $\G_1$ contains a bad vertex $v$.
Moreover, since $\G_1$ is connected by Proposition \ref{handa}(iii), we can further assume that 
$v \in N(w)$ for some good vertex $w \in \G_1 \cup \{a,b\}$.
Note that since $\G$ is bipartite, we have that $d(v,a) > d(w,a)$ and $d(v,b) > d(w,b)$. 
Let $y \in \G_2$. It follows from Lemma \ref{lem1} that $d(w,y) < d(v,y)$,
and so $y \in W_{wv}$. Therefore $W_{vw} \subseteq \G_1$, implying that
$$
  \frac{n}{2} = |W_{vw}| \le |\G_1|.
$$
\qed

\begin{corollary}
\label{cor3}
Let $\G$ be as in Convention \ref{blank}.
Then either $\G_1$ or $\G_2$ contains no bad vertex.
\end{corollary}
\proof
Suppose that both $\G_1$ and $\G_2$ contain bad vertices. 
By Lemma \ref{lem2} we get
$$
  n-2 = |\G_1|+|\G_2| \ge \frac{n}{2} + \frac{n}{2} = n,
$$
a contradiction. \qed

\begin{lemma}
\label{lem4}
Let $\G$ be as in Convention \ref{blank}.
Then either $\G_1$ or $\G_2$ contains no good vertex.
\end{lemma}
\proof
Corollary \ref{cor3} shows that one of $\G_1$ and $\G_2$ consists entirely of 
good vertices. Without loss of generality we can assume that this holds
for $\G_2$.
Now, suppose that $x' \in \G_1$ is a good vertex.
Then there also exists a good vertex $x \in N(b) \cap \G_1$.
Let $Z$ be the set of vertices $z$ of $\G_1$, for which at least one shortest $z$-$a$ path of 
the subgraph $\la \G_1 \cup \{a,b\} \ra$ passes through $b$. Observe that 
$\G_2 \cup Z \cup \{b\} \subseteq W_{bx}$. 
Let $B_1 = N(a) \cap \G_2$ and pick $y \in B_1$. Since $x$ and $y$ are both good vertices, we have that $d(a,x) = d(b,y)$, and so
$W_{ya} \subseteq \G_2 \setminus B_1 \cup Z \cup \{b,y\}$. Hence
$$
  |\G_2| + |Z| + 1 \le \frac{n}{2} \le |\G_2| - |B_1| + |Z| + 2,
$$
implying that $|B_1|=1$. This contradicts the minimality of $d(a,b)$. \qed

\begin{figure}[h]
\begin{center}
\includegraphics[scale=0.6]{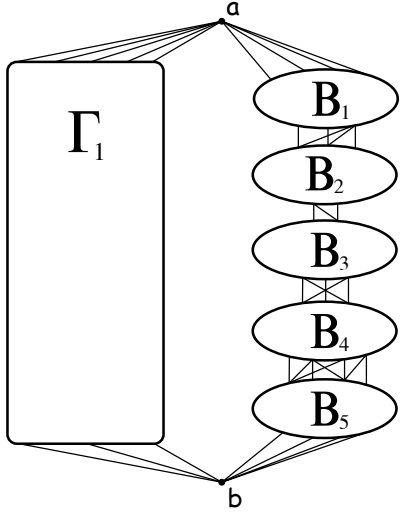}
\caption{The sets $B_i$ of the good component.}
\label{fig:Bji}
\end{center}
\end{figure}

The above corollary and lemma enable us to make the following convention
for the rest of the paper.
\begin{convention}
\label{blank1}
{\em Let $\G$ be as in Convention \ref{blank}.
By Corollary \ref{cor3} and Lemma \ref{lem4} one of $\G_1, \G_2$ coincides
with the set of all bad vertices. Without loss of generality we assume that this holds for $\G_1$.
Consequently, $\G_2 \cup \{a,b\}$ coincides with the set of all good vertices.
We call $\G_1$ the {\em bad component} and $\G_2$ the {\em good component} of $\G$.
Furthermore, Corollary~\ref{cor3} implies that the vertex set of $\G_2$ can be partitioned into $m$ sets $B_1 \cup B_2 \cup \cdots \cup B_m$, where $m = d(a,b) - 1$ and 
$B_i = N_i(a) \cap \G_2 = N_{m+1-i}(b) \cap \G_2$. Moreover, for every $x \in B_i$, where $1 \leq i \leq m$, we have that $N(x) \subseteq B_{i-1} \cup B_{i+1}$ 
(with the agreement that $B_0 = \{a\}$ and $B_{m+1} = \{b\}$).
See Figure~\ref{fig:Bji} for the case $m = 5$.}
\end{convention}

\section{The main results}
\label{sec:goodcomp}

Let $\G$ be as in Convention \ref{blank1}. In this section we show that $d(a,b) \ge 3$ and give a complete
classification of the graphs $\G$ for which $d(a,b)=3$.
We first show that every $x \in B_1$ and every $y \in B_m$ lie on a common shortest $a$-$b$ path of $\G$.

\begin{lemma}
\label{lem6}
Let $\G$ be as in Convention \ref{blank1}.
Then for each $x \in B_1$ and $y \in B_m$ we have $d(x,y) = d(a,b) - 2$.
\end{lemma}
\proof
Let $T_1 = W_{xa} \cap \G_1$ and let $T_2 = W_{xa} \cap \G_2$. Since $W_{xa} = T_1 \cup T_2 \cup \{b\}$, it thus follows that $\frac{n}{2} = |T_1| + |T_2| + 1$. Note that, 
by Corollary~\ref{cor3}, it suffices to prove that $y \in T_2$.

Suppose to the contrary that $y \notin T_2$. We show that in this case $T_2 \subseteq W_{by}$. Indeed, let $z \in T_2$ and let $i$ be such that $z \in B_i$.
Then Corollary~\ref{cor3} implies that $d(z,b) = m-i+1$. On the other hand $d(z,y) = m-i+2$, since otherwise $y \in T_2$. 
Of course, $T_1 \cup \{b\} \subseteq W_{by}$, and so  $\frac{n}{2} = |T_1| + |T_2| + 1$ implies that 
\begin{equation}
\label{eq:lem6}
	W_{by} = T_1 \cup T_2 \cup \{b\}.
\end{equation} 
Now, since $a$ has at least one neighbor in $\G_1$ (which, of course, does not belong to $T_1$), 
there exists at least one vertex $w \in \G_1 \setminus T_1$, adjacent to some vertex $v \in T_1$. 
By \eqref{eq:lem6} we have that $w \in W_{yb}$, that is $d(w,y) = d(w,a) + d(a,y)$, and $d(w,b) = d(w,y)+1$. 
However, since $v \in T_1$, Proposition \ref{handa}(iv) implies that 
$$
	\begin{array}{lcl}
	d(w,y) & = & d(w,a) + d(a,y) = d(v,x) + d(a,y) = d(v,b) + d(b,x) + d(a,y) \\
		&  = & d(w,b) - 1 + d(b,x) + d(a,y) \geq d(w,b) + 1,
  \end{array}
$$
contradicting $w \in W_{yb}$. Hence, $y \in T_2$ as claimed. \qed

\begin{remark}
\label{rem:d=3}
{\em 
Let $\G$ be as in Convention \ref{blank1}.
Note that Lemma \ref{lem6} implies that, in the case of $d(a,b)=3$, we have that $\la B_1 \cup B_2 \ra$ is the complete bipartite graph $K_{|B_1|,|B_2|}$.}
\end{remark}

\begin{proposition}
\label{pro:razdalja3}
Let $\G$ be as in Convention \ref{blank1}.
Then $d(a,b) \geq 3$.
\end{proposition} 
\proof
By Proposition \ref{handa}(i), $d(a,b) \ge 2$.
Observe that, if $d(a,b) = 2$, then Lemma~\ref{lem6} implies that $\G_2$ consists of a single vertex 
(note that this also follows by Corollary~\ref{cor3} and Proposition \ref{handa}(iii)). 
But then the minimal degree of $\G$ is less than $3$, contradicting Proposition \ref{handa}(ii). \qed

\noindent
For the rest of this section
we introduce the following notation. Let $\G$ be as in Convention \ref{blank1}.
Assume that in addition every vertex of $\G_1$ lies on some shortest $a$-$b$ path 
of the subgraph $\wtg=\la \G_1 \cup \{a,b\} \ra$.
Then the vertex set of $\wtg$ can be partitioned into $t+2$ sets $D_0 \cup D_1 \cup \cdots \cup D_{t+1}$, where $t + 1$ is 
the distance between $a$ and $b$ in $\wtg$, and $D_i = \widetilde{N}_i(a) = \widetilde{N}_{t+1-i}(b)$ 
(where by $\widetilde{N}_i(a)$ (resp. $\widetilde{N}_i(b)$) we mean the $i$-th
neighbourhood of $a$ (resp. $b$) in $\wtg)$. Moreover, for every $x \in D_i$, where $1 \leq i \leq t$, 
we have that $N(x) \subseteq D_{i-1} \cup D_{i+1}$.

\begin{lemma}
\label{lem:d=3}
Let $\G$ be as in the above paragraph and assume in addition that $d(a,b)=3$. 
Then the number $t$ is even.
Let $Y = D_1 \cup \cdots \cup D_{t/2-1}$ and $X = D_{t/2+2} \cup \cdots \cup D_t$.
Then the following hold:
\begin{itemize}
\item[(i)]
$|D_{t/2}| = |D_{t/2+1}| = 1$.
\item[(ii)]
$|B_1|=|B_2|$ and $|X| = |Y|$.
\item[(iii)]
$|D_{t/2-1}| = |D_{t/2+2}| = |B_1|$, $|D_1|=|D_t|=1$ and $t \geq 8$. 
\item[(iv)]
Every vertex of $D_{t/2+2}$ is adjacent to every vertex of $D_{t/2+3}$
and every vertex of $D_{t/2-1}$ is adjacent to every vertex of $D_{t/2-2}$.
\item[(v)]
$|D_{t/2+4}| = |D_{t/2-3}| = 1$.
\end{itemize}
\end{lemma}
\proof
That $t$ is even is clear as $\G$ is bipartite. Since $\G_1$ contains no good vertices, this immediately implies $t \ge 4$. \\
(i) Pick $x \in B_1$, $y \in B_2$. 
Recall that, by Remark \ref{rem:d=3}, $\la B_1 \cup B_2 \ra$ is a complete bipartite graph. 
Note that $X = W_{xa} \cap \G_1$ and $Y = W_{yb} \cap \G_1$. Let $Z = \G_1 \setminus (X \cup Y)$, that is,  
$Z = D_{t/2} \cup D_{t/2+1}$. Since $W_{xa} = X \cup B_2 \cup \{b,x\}$ and $W_{yb} = Y \cup B_1 \cup \{a,y\}$, 
we have that $W_{xa} \cup  W_{yb} = V(\G) \setminus Z$. As $W_{xa} \cap W_{yb} = \{x,y\}$, we thus have that 
$n - |Z| = |W_{xa} \cup W_{yb}| = \frac{n}{2} + \frac{n}{2} - 2$, and so $|Z| = 2$.
Since $D_i \ne \emptyset$ by definition, the result follows.

\smallskip \noindent
(ii) Let $\tilde{y}$ be the unique vertex in $D_{t/2}$ and let $\tilde{x}$ be the unique vertex in $D_{t/2+1}$.
Note that $W_{\tilde{y} \tilde{x}} = \{\tilde{y},a\} \cup Y \cup B_1$. The result follows since 
$W_{xa} = X \cup B_2 \cup \{b,x\}$ and $W_{xy} = B_2 \setminus \{y\} \cup \{a,x,\tilde{y}\} \cup Y$.

\smallskip \noindent
(iii) Pick $v \in D_{t/2+2}$ and observe that $D_{t/2+2} \setminus \{v\} \cup Y \cup \{\tilde{x},\tilde{y},a\} \subseteq W_{\tilde{x} v}$.
Therefore, $|D_{t/2+2}| + |Y| + 2 \le n/2$. Since $n/2 = |W_{\tilde{y} \tilde{x}}| = 2 + |Y| + |B_1|$, this implies $|D_{t/2+2}| \le |B_1|$.
Pick now $w \in D_1$ and observe that $D_1 \setminus \{w\} \cup \{a,b\} \cup B_1 \cup B_2 \cup X \setminus D_{t/2+2} \subseteq W_{a w}$.
Therefore $|D_1|-1 + 2 + |B_1| + |B_2| + |X| - |D_{t/2+2}| \le n/2$.
As $n/2 = |W_{\tilde{x} \tilde{y}}| = 2 + |X| + |B_2|$, this implies $|D_1| - 1 + |B_1| - |D_{t/2+2}| \le 0$.
But since $|D_{t/2+2}| \le |B_1|$, we must have $|D_1| = 1$ and $|B_1| = |D_{t/2+2}|$.
Similarly we find that $|D_t| = 1$ and $|D_{t/2-1}| = |B_2| = |B_1|$.
Observe that this immediately implies $t \ge 6$ since otherwise $D_1 = D_{t/2-1}$, and so $|B_1|=1$, contradicting minimality of $d(a,b)$.
Moreover, as $|D_1|=1$, $t\ne 6$, since otherwise the distance between $\tilde{y}$ and the unique vertex of $D_1$ is smaller than $d(a,b)$.
Hence, $t \ge 8$ as claimed.

\smallskip \noindent
(iv) Pick $v \in D_{t/2+2}$ and suppose there is $w \in D_{t/2+3}$ which is not adjacent to $v$.
Observe that in this case we have $D_{t/2+2} \setminus \{v\} \cup \{w, \tilde{x}, \tilde{y}, a\} \cup Y \subseteq W_{\tilde{x} v}$.
As $|D_{t/2+2}| = |B_1|$ by (iii) above this implies $3 + |B_1| + |Y| \le n/2$.
But $n/2 = |W_{\tilde{y} \tilde{x}}| = 2 + |Y| + |B_1|$, a contradiction. Hence $v$ is adjacent to all vertices of $D_{t/2+3}$.
Similarly we show that every vertex of $D_{t/2-1}$ is adjacent to every vertex of $D_{t/2-2}$.

\smallskip \noindent
(v) The case $t=8$ is covered by (iii) above. We can thus assume that $t \ge 10$. 
Pick adjacent vertices $v \in D_{t/2+3}$ and $w \in D_{t/2+4}$.
Observe that, by (iv) above, we have 
$$
  (N(v) \cap D_{t/2+4}) \setminus \{w\} \cup \{v, \tilde{x},\tilde{y}\} \cup D_{t/2+2} \cup (Y \setminus D_1) \cup (D_{t/2+3} \setminus N(w)) \subseteq W_{vw}.
$$
Therefore 
$$
  |N(v) \cap D_{t/2+4}| + 2 + |D_{t/2+2}| + |Y| - 1 + |D_{t/2+3} \setminus N(w)| \le n/2.
$$
Combining together (iii) above and $2 + |Y| + |B_1| = n/2$, we find that
$|N(v) \cap D_{t/2+4}| - 1 + |D_{t/2+3} \setminus N(w)| \le 0$. Consequently, $|N(v) \cap D_{t/2+4}| = 1$ and $|D_{t/2+3} \setminus N(w)|=0$.
In particular, $D_{t/2+3}\subseteq N(w)$. The result follows. Similarly we show that $|D_{t/2-3}| = 1$. \qed

\begin{figure}[h]
\begin{center}
\includegraphics[scale=0.8]{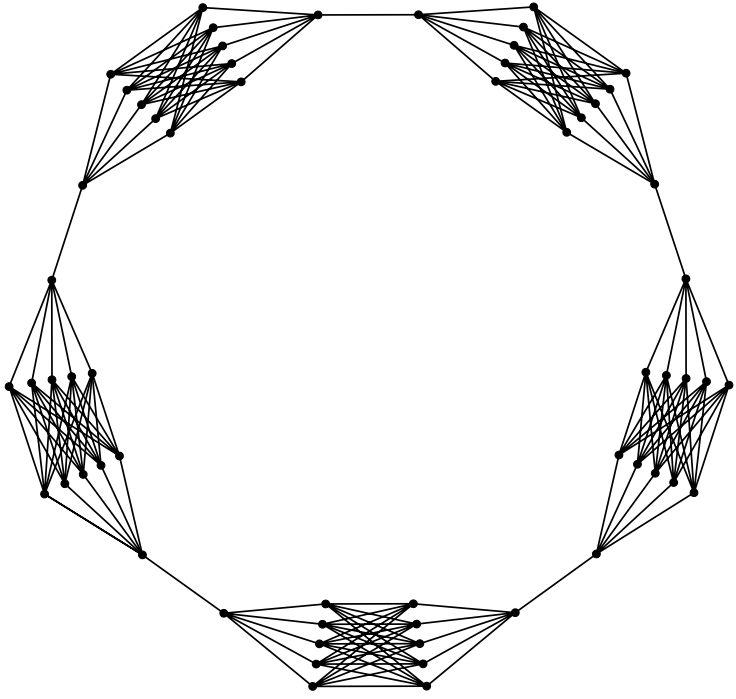}
\caption{The graph $W(5,5)$.}
\label{fig:Wheel-sun}
\end{center}
\end{figure}

Before stating the classification of graphs satisfying the assumptions of Convention~\ref{blank1} with $d(a,b)=3$ we give the following construction.
\begin{construction}
{\em Let $m \ge 2$ and $\ell \ge 3$ be integers. For $0 \le i \le 4 \ell-1$ let $f(i)$ be $1$ if $i$ is congruent to $0$ or $3$ modulo $4$, and $m$ otherwise.
The graph $W(m, \ell)$ has vertex set 
$$
  V(W(m,\ell)) = \{ (i,j) \mid 0 \le i \le 4\ell-1, 1 \le j \le f(i) \}
$$
and edge set
$$
  E(W(m,\ell)) = \{ \{(i_1,j_1), (i_2, j_2)\} \mid i_2-i_1 \equiv 1 (\hspace{-3mm}\mod 4\ell)\}.
$$
In other words, the graph $W(m,\ell)$ is obtained from the cycle of length $4 \ell$ by replacing every second pair of vertices by
a complete bipartite graph $K_{m,m}$, see Figure \ref{fig:Wheel-sun}.}
\end{construction}

\begin{theorem}
\label{the:Wheel-sun}
Let $m \ge 2$ and $\ell \ge 3$ be integers. The graph $W(m,\ell)$ is bipartite and not $3$-connected. Moreover, it is distance-balanced if and only if $\ell$ is odd.
\end{theorem}
\proof
That $W(m,\ell)$ is bipartite and not $3$-connected is clear from the construction. To show that $W(m,\ell)$ is not distance-balanced in the case when $\ell$ is even one only needs to observe that the edges of the form $\{(4i,j_1),(4i+1,j_2)\}$ and $\{(4i+2,j_1),(4i+3,j_2)\}$ are ``non-balanced''. To show that $W(m,\ell)$ is distance-balanced when $\ell$ is odd first observe that the automorphism group of $W(m,\ell)$ has three orbits on the edge set, the orbit of $\{(0,1),(1,1)\}$, the orbit of $\{(1,1),(2,1)\}$ and the orbit of $\{(3,1),(4,1)\}$. To check that each of these three edges is ``balanced'' is an easy exercise.
\qed

\begin{theorem}
\label{thm:main}
Let $\G$ be a bipartite distance-balanced graph for which there exists a $2$-cut $\{a,b\}$ with $d(a,b)=3$.
Then either $\G$ is a cycle or there exists an integer $m \ge 2$ and an odd integer $\ell \ge 3$ such that $\G$ is isomorphic to the graph $W(m,\ell)$.
\end{theorem} 
\proof
Assume that $\G$ is not a cycle and note that in this case $\G$ satisfies the assumptions of Convention \ref{blank1}.
Pick $x \in B_1$, $y \in B_2$. 
By Lemma~\ref{lem6} we have that $B_1 = N(y) \setminus {b}$ and $B_2 = N(x) \setminus \{a\}$. 
Let $X = W_{xa} \cap \G_1$, let $Y = W_{yb} \cap \G_1$ and let $Z = \G_1 \setminus (X \cup Y)$. 
Observe that for any $w \in X$ every shortest $w$-$x$ path of $\G$ passes through $b$. Hence $w \in W_{by}$, 
implying that $X \cap Y = \emptyset$. Moreover, no vertex of $X$ is adjacent to some vertex of $Y$.
Namely, if $w \in X$ is adjacent to $v \in Y$ then $w \in W_{yx}$ and $v \in W_{xy}$, and so Proposition \ref{handa}(iv)
implies
$$
  d(w,b) + 1 = d(w,y) = d(v,x) = d(v,b) - 2 \le d(w,b)-1,
$$
a contradiction. Since 
\begin{equation}
\label{eq:d=3b}
W_{xa} = X \cup B_2 \cup \{b,x\} \qquad \hbox{and} \qquad W_{yb} = Y \cup B_1 \cup \{a,y\}, 
\end{equation}
we have that $W_{xa} \cup  W_{yb} = V(\G) \setminus Z$. As $W_{xa} \cap W_{yb} = \{x,y\}$, we thus have that 
$$
 	n - |Z| = |W_{xa} \cup W_{yb}| = \frac{n}{2} + \frac{n}{2} - 2,
$$
and so $|Z| = 2$.

\medskip \noindent
{\bf Claim 1:} For $z \in Z$ we have that either $N(z) \cap X = \emptyset$ or $N(z) \cap Y = \emptyset$.

\noindent
Let $z \in Z$ be such that $N(z) \cap Y \ne \emptyset$ and
pick $v \in N(z) \cap Y$. Since $v \in W_{yb}$ and $z \in W_{by}$, Proposition \ref{handa}(iv) implies that $d(z,b) = d(v,y)$. 
Since all shortest $v$-$y$ paths pass through $a$ and $\G$ is bipartite we thus have
$d(z,b) = d(v,y) = d(v,a) + d(a,y) = d(z,a)+1$.
Similarly, if $N(z) \cap X \ne \emptyset$, we find that  
$d(z,a) = d(z,b)+1$. This proves Claim 1, since otherwise $d(z,a) = d(z,a)+2$.

\medskip \noindent
We may therefore assume that $Z=\{\tilde{x}, \tilde{y}\}$ where $N(\tilde{y}) \cap X = \emptyset$ and $N(\tilde{x}) \cap Y = \emptyset$.
Since $\G_1$ is connected $\tilde{x}$ and $\tilde{y}$ are adjacent and $N(\tilde{x}) \cap X \ne \emptyset$ and $N(\tilde{y}) \cap Y \ne \emptyset$.
The arguments above (for $z = \tilde{y}$ and $z = \tilde{x}$) yield $d(\tilde{y},b) = d(\tilde{y},a)+1$ and $d(\tilde{x},a) = d(\tilde{x},b)+1$.
Since $\tilde{y} \in W_{by}$ and $N(\tilde{y}) \cap X = \emptyset$, we have that $d(\tilde{y},b) = d(\tilde{x},b)+1$, and so 
\begin{equation}
\label{eq:d=3}
  d(\tilde{x},b)=d(\tilde{y},a).
\end{equation}
Observe also that, since $\G \setminus \{a,\tilde{y}\}$ and $\G \setminus \{b,\tilde{x}\}$ are disconnected, we have 
$d(\tilde{x},b)=d(\tilde{y},a) \ge 3$ by minimality of $d(a,b)$.

\medskip \noindent
{\bf Claim 2:} every vertex of $\G_1$ lies on some shortest $a$-$b$ path 
in the subgraph $\wtg= \la \G_1 \cup \{a,b\} \ra$.\\
Suppose that there is a vertex $u \in X$ which does not lie on such a path. Since $\G_1$ is connected we may assume that
there is a neighbour $w$ of $u$ which does lie on such a path.
But since $\G$ is bipartite we now have $W_{uw} \subseteq X \setminus \{w\}$, implying that $n/2 = |W_{uw}| \le |X|-1 < |X|$,
which, by \eqref{eq:d=3b}, is impossible.
Similarly we show that every vertex in $Y$ lies on some shortest $a$-$b$ path of $\wtg$. This proves Claim 2.

\smallskip \noindent
Let $t$ and sets $D_i \; (0 \le i \le t+1)$ be as in the paragraph preceding Lemma~\ref{lem:d=3}. 
By Lemma~\ref{lem:d=3} we have that $t \ge 8$. Let $v$ be the unique vertex of $D_{t/2+4}$ and
observe that the pair $(\tilde{x},v)$ satisfies the same assumptions as the pair $(a,b)$. We can thus apply 
Lemma \ref{lem:d=3} to it. Hence 
we find that $|D_{t/2+3}| = |D_{t/2+2}| = |B_1| = |D_2|$, $|D_4| = |D_{t/2+5}| = 1$, and that every vertex of 
$D_2$ is adjacent to every vertex of $D_3$. 
Note that as $|D_4|=1$, also $|D_{t/2+8}|=1$ (again by replacing the pair $(a,b)$ by $(\tilde{x},v)$).
Continuing with this process we find that the following two possibilities can occur:
\begin{itemize}
\item[(i)]
There is some $i$ such that $|D_i| = |D_{i+1}| = |D_{i+2}| = 1$
(note that $i$ can be $0$). But this contradicts minimality of $d(a,b)$.
\item[(ii)]
$t$ is divisible by $8$ and $|D_{2+4j}| = |D_{3+4j}| = |B_1| = |B_2|$ for $0 \le j \le t/4-1$,
while $|D_j| = 1$ for all other $j$. Therefore $n = 2|B_1| (t/4+1) + 2(t/4+1) = 2(t/4+1)(|B_1|+1)$.
It is now clear that $\G$ is isomorphic to the graph $W(m,\ell)$, where $m=|B_1|$ and $\ell = t/4+1$.
\end{itemize}
\qed

\begin{figure}[h]
\begin{center}
\includegraphics[scale=1]{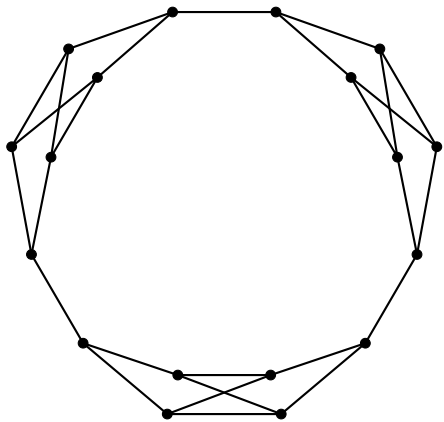}
\caption{The smallest bipartite non-$3$-connected distance-balanced graph, which is not a cycle.}
\label{fig:smallest Wheel-sun}
\end{center}
\end{figure}

\begin{theorem}
\label{the:najmanjsi}
Except for cycles of even length, the order of a bipartite non-$3$-connected distance-balanced graph
of smallest order is $18$. Moreover, such a graph is unique and  
is isomorphic to the graph $W(2,3)$ (see Figure~\ref{fig:smallest Wheel-sun}).
\end{theorem}
\proof
Let $\G$ be a bipartite non-$3$-connected distance-balanced graph, which is not a cycle, of smallest order. 
By  Proposition~\ref{pro:razdalja3} we have that $d(a,b) \geq 3$, where we use the notation from Convention~\ref{blank1}. 
Theorem~\ref{the:Wheel-sun} implies that $n \le 18$. 
If $d(a,b) = 3$ then Theorem~\ref{thm:main} applies. We can thus assume that $d(a,b) \geq 4$. 
Let the sets $D_i$, $1 \leq i \leq t$, represent the vertices of $\G_1$ at distance $i$ from $a$ 
in the subgraph $\wtg=\la \G_1 \cup \{a,b\} \ra$, where $t$ is the maximal distance from $a$ in $\wtg$. 
Since $\G_1$ contains no good vertex, we have that the distance between $a$ and $b$ in $\wtg$ is at least $m+3$, 
and so $t \geq m+2$ (recall that $d(a,b) = m+1$). Moreover, by minimality of $d(a,b)$ we have that $|D_i| \geq 2$ 
for all $2 \leq i \leq m$ and that $|B_i| \geq 2$ for all $1 \leq i \leq m$. Therefore, $\G$ is of order at least $4m+3$. 
Thus, if $d(a,b) \geq 5$, $\G$ is of order at least $19$. We are therefore left with the case $d(a,b) = 4$. 
Pick $x \in B_1$ and observe that $B_1 \setminus \{x\} \cup D_1 \cup D_2 \cup D_3 \cup D_4 \cup \{a\} \subseteq W_{ax}$. 
Since $D_1$ and $D_4$ cannot both be singletons (otherwise $d(a,b)$ is not minimal), this implies 
$\frac{n}{2} = |W_{ax}| \geq 9$. Hence, $n = 18$ and $W_{ax} = B_1 \setminus \{x\} \cup D_1 \cup D_2 \cup D_3 \cup D_4 \cup \{a\}$. Moreover, 
$|B_1| = |D_2| = |D_3| = 2$ and $|D_1 \cup D_4| = 3$. A similar argument (considering the distance partition from $b$ in $\wtg$) shows that also $|B_3| = 2$.
Note that this implies that $|D_5 \cup D_6 \cup \cdots \cup D_t| \le 3$, and so Proposition \ref{handa}(ii) implies $t \le 6$ 
(otherwise the degree of the unique vertex of $D_6$ is $2$). If $|D_4|=1$, say $D_4=\{c\}$, then $\{c,b\}$ is a $2$-cut with $d(b,c) \le 3$, a contradiction.
It follows that $|D_4|=2$ and $|D_1|=1$.
Let now $w$ be the unique vertex of $D_1$ and let $v \in D_2$. Then $D_2 \setminus \{v\} \cup \{w,a\} \cup B_1 \cup B_2 \cup B_3 \subseteq W_{wv}$, and so
$|W_{wv}| \geq 9$, implying that $|B_2| = 2$. Since the minimal degree in $\G$ is at least $3$, every vertex of $D_2$ is adjacent to every vertex of $D_3$. 
We show that $b$ has no neighbour in $D_5$. Indeed, if $b$ has a neighbor in $D_5$, every vertex of $D_2$ is closer to $b$ than to any vertex of $B_3$. 
Hence, for $u \in B_3$ we have that $W_{ub} = B_2 \cup B_1 \cup \{u,a,w\}$, and so $|W_{ub}| = 7$, a contradiction. It follows that $t = 7$ and $|D_5| = |D_6| = |D_7| = 1$, 
implying that the unique vertex of $D_6$ is of degree $2$, a contradiction. This completes the proof.
\qed

\section{Further structural properties}
\label{sec:badcomp}

Throughout this section let $\G$ be as in Convention~\ref{blank1}. 
To further describe the structure of the bad component we divide the edges of $\wtg = \la \G_1 \cup \{a,b\}\ra$ into two classes.
Let $\wtd$ denote the distance function of the subgraph $\wtg$.
The edge $xy$ of $\wtg$ for which either $\wtd (x,a) = \wtd (y,a) + 1$ and $\wtd (x,b) = \wtd (y,b) + 1$, or
$\wtd (y,a) = \wtd (x,a) + 1$ and $\wtd (y,b) = \wtd (x,b) + 1$, is called {\em horizontal}.
In this case the vertex of $xy$ that is closer to $a$ (and also to $b$) 
is called the {\em right} vertex of $xy$. The other vertex of $xy$ is called the {\em left} vertex of $xy$. 
All other edges of $\wtg$ are called {\em vertical}. Note that, since $\G$ is bipartite, for a vertical edge $xy$ we either have
$\wtd (x,a)+1 = \wtd (y,a)$ and $\wtd (x,b) = \wtd (y,b) + 1$, or $\wtd (y,a)+1 = \wtd (x,a)$ and $\wtd (y,b) = \wtd (x,b) + 1$. 
In this case the vertex of $xy$ that is closer to $a$ (and further from $b$) is called the {\em upper} vertex of $xy$. The other vertex of $xy$ is called the {\em lower} vertex of $xy$. 

\begin{proposition}
\label{prop:horvert}
Let $\G$ be as in Convention~\ref{blank1}. Then the following (i), (ii) hold.
\begin{itemize}
\item[(i)]
Let $xy$ be a horizontal edge, where $x$ is the left and $y$ the right vertex of $xy$.
Then $W_{xy} \subseteq \G_1$.
\item[(ii)]
Let $xy$ be a vertical edge, where $x$ is the upper and $y$ the lower vertex of $xy$.
Let 
\begin{equation}
\label{eq:ell}
  \ell = (\wtd(a,b) + \wtd(x,a) - \wtd(y,b) + 1)/2.
\end{equation}
Then
\begin{equation}
\label{formula}
  |W_{xy} \cap \G_1| = \frac{n}{2} - \sum_{j=0}^{(\wtd(a,b) + m + 1)/2 - \ell} |B_j|,
\end{equation}
with the agreement that $B_j = \emptyset$ for $j < 0$ and $j > m+1$.   
\end{itemize}
\end{proposition}
\proof
(i) Since $y$ is closer to both $a$ and $b$ than $x$, Lemma \ref{lem1} implies that $\G_2 \cup \{a,b\} \subseteq W_{yx}$.
Therefore $W_{xy} \subseteq \G_1$.

\smallskip \noindent
(ii) Consider the set $B_j \; (0 \le j \le m+1)$. By Corollary~\ref{cor3} all the vertices of this set
have the same distance to $a$, and all have the same distance to $b$. Consequently, either $B_j \subseteq W_{xy}$ or
$B_j \cap W_{xy} = \emptyset$. Set $j = (\wtd(a,b) + m + 1)/2 - \ell$.

Suppose first that $j \ge m+1$, that is, $2 \ell \le \wtd(a,b) - m - 1$.
By \eqref{eq:ell} we have that $\wtd(x,a) + m + 1 \le \wtd(y,b) - 1$, and so $b$, and hence also $\cup_{j=0}^{m+1} B_j$,
is contained in $W_{xy}$. The equation \eqref{formula} thus holds in this case.

Suppose next that $j \le -1$, that is $\wtd(a,b)+m+3 \le 2 \ell$.
By \eqref{eq:ell} we have that $\wtd(b,y) + m + 1 \le \wtd(x,a)-1$, and so $a$, and hence also $\cup_{j=0}^{m+1} B_j$,
is contained in $W_{yx}$. The equation \eqref{formula} thus holds in this case.

Suppose finally that $0 \le j \le m$. 
Pick $u \in B_j$ and $v \in B_{j+1}$.
Similarly as above we show that $u \in W_{xy}$, while $v \in W_{yx}$.
It follows that $B_i \subseteq W_{xy}$ for $0 \le i \le j$ and 
$B_i \subseteq W_{yx}$ for $j+1 \le i \le m+1$.
This completes the proof. \qed

\section{Applications}
\label{sec:appl}

In this section we give two applications of the results from
the previous sections. In particular, we show that if we restrict either to the family of bipartite strongly distance-balanced graphs or to the family of distance-balanced partial cubes, the answer to the Handa question \ref{que} is affirmative. Note that the result for distance-balanced partial cubes was
already obtained by Handa~\cite{Ha} - here we give an alternative proof.

We first consider the strongly distance-balanced graphs. 
A graph $\G$ is {\em strongly distance-balanced} if and only if for every pair of adjacent vertices $u,v \in V(\G)$ 
and for any nonnegative integer $i$ we have $|N_i(u) \cap N_{i+1}(v)| = |N_{i+1}(u) \cap N_i(v)|$.
Clearly, every strongly distance-balanced graph is distance-balanced.
Strongly distance-balanced graphs were first introduced in~\cite{KMMM}. It turns out (see~\cite{KMMM}) that strongly distance-balanced graphs actually coincide with the so-called distance degree regular graphs (introduced in~\cite{HN84}). In~\cite{KMMM} it was also observed 
that not every distance-balanced graph is strongly distance-balanced. (As a corollary of Theorem~\ref{the:Wheel-sun} and the next theorem we have that every graph $W(m,\ell)$, where $\ell \geq 3$ is odd, is a distance-balanced graph that is not strongly distance-balanced.)

\begin{theorem}
\label{thm:sdb}
Let $\G$ be a bipartite strongly distance-balanced graph that is not a cycle. Then $\G$ is $3$-connected.
\end{theorem}
\proof
Suppose to the contrary that $\G$ is not $3$-connected. Then $\G$ satisfies the assumptions of Convention~\ref{blank1}.
Pick $x \in B_1$ and consider the set $N_{m+1}(a) \cap N_m(x)$.
By Corollary \ref{cor3} we have $N_{m+1}(a) \cap N_m(x) = \{b\}$. 
Therefore, since $\G$ is strongly distance-balanced, $|N_m(a) \cap N_{m+1}(x)| = 1$. 
By Lemma \ref{lem6}, $N_m(a) \cap N_{m+1}(x) \subseteq \G_1$.
However, every vertex of $\G_1$, which is at distance $m$ from $a$ is at distance $m+1$ from $x$.
Hence $N_m(a) \cap N_{m+1}(x) = N_m(a) \cap \G_1$, and so $|N_m(a) \cap \G_1| = 1$. Let $c$ be the unique vertex of $N_m(a) \cap \G_1$.
Observe that $\G \setminus \{a,c\}$ is disconnected. For example, every path from $b$ to some 
vertex in $N(a) \cap \G_1$ passes either through $a$ or through $c$.
But this contradicts the minimality of $d(a,b)$. \qed

A graph $\G$ is a {\em partial cube} if it can be isometrically embedded into a hypercube.
It was proven in \cite{Dj} that $\G$ is a partial cube if and only if it is 
bipartite and for every pair of adjacent vertices $u$ and $v$ of $\G$, the set $W_{uv}$ is convex 
(the set $U \subseteq V(\G)$ is {\em convex} if for every $x,y \in U$ every shortest  
$x$-$y$ path is contained in $U$). It follows from \cite[Theorem 1.3]{Ha} that if $\G$ is a distance-balanced
partial cube that is not a cycle, then $\G$ is $3$-connected. We now give an alternative proof of this result.
 
\begin{theorem}
\label{thm:pc}
Let $\G$ be a distance-balanced partial cube that is not a cycle. Then $\G$ is $3$-connected.
\end{theorem}
\proof
Suppose to the contrary that $\G$ is not $3$-connected. Then $\G$ satisfies the assumptions of Convention~\ref{blank1}. Note that, by the minimality of 
$d(a,b)$, we have $|B_1|\ge 2$. Pick $x,y \in B_1, x \ne y$, and consider the set
$W_{ax}$. Clearly $y \in W_{ax}$. Note that since all vertices of $\G_1$ are bad vertices,   
$W_{xa} \cap \G_1 \ne \emptyset$. Pick adjacent vertices $v,w \in \G_1$ such that $v \in W_{ax}$ and $w \in W_{xa}$
(such vertices exist since $\G_1$ is connected). By Lemma \ref{lem1} we have
$d(v,y) = \min \{d(v,a) + 1, d(v,b) + d(b,y)\}$. But by Proposition \ref{handa}(iv) we have
$d(v,a) + 1 = d(w,x) + 1 = d(w,b) + d(b,x) + 1 = d(w,b) + d(b,y) + 1 = d(v,b) + d(b,y)$, implying that there exists a shortest
$y$-$v$ path passing through $b$. Since $\G$ is a partial cube and $y,v \in W_{ax}$,
this implies $b \in W_{ax}$, contradicting Corollary~\ref{cor3}. \qed

\section{Directions for future research}
\label{sec:research}

Although the Handa question has now been answered, there are still numerous possibilities for further research.
In this section we propose some of the questions and problems that we believe are worth considering.
We continue to use the notation from Convention~\ref{blank1}.

In view of the fact that Theorem~\ref{the:Wheel-sun} and Theorem~\ref{thm:main} give a complete classification of 
the graphs with $d(a,b)=3$, the first natural question is the following.

\begin{question}
Are there any bipartite, non-$3$-connected distance-balanced graphs, which are not cycles, for which
the minimal distance between the vertices of a $2$-cut is at least $4$? 
\end{question}

In the case that the answer to the above question is positive, the following problem should be considered.

\begin{problem} 
Classify the bipartite, non-$3$-connected distance-balanced graphs.
\end{problem}

If this problem turns out to be to difficult to solve in general, one could restrict to particular families of 
bipartite graphs (such as in Theorems \ref{thm:sdb} and \ref{thm:pc}).

The next direction one should try to pursue is the investigation of the connectivity of non-bipartite
distance-balanced graphs (recall that, by Proposition~\ref{handa}, such graphs are at least $2$-connected).
Of course, a complete classification of all non-$3$-connected distance-balanced graphs seems out of reach in
view of the above problem. One should thus first try to generalize the results of \cite{Ha} and this paper to
non-bipartite graphs.



\begin{thebibliography}{99}

\begin{footnotesize}


\bibitem{BCPSSS} K.~Balakrishnan, M.~Changat, I.~Peterin, S.~\v{S}pacapan, P.~\v{S}parl, A.~R.~Subhamathi, 
                 Strongly distance-balanced graphs and graph products, 
                 {\em European J. Combin.} {\bf 30} (2009), 1048--1053. 
                 
\bibitem{Dj}     D.~\v{Z}.~Djokovi\'c, Distance preserving subgraphs of hypercubes, 
                 {\em J. Combin. Theory Ser. B} {\bf 14} (1973), 263--267.

\bibitem{Ha} K.\ Handa, Bipartite graphs with balanced $(a,b)$-partitions,
                {\em Ars Combin.} {\bf 51} (1999), 113--119.
                
\bibitem{HN84} T.~Hilado, K.~Nomura, Distance Degree Regular Graphs, 
               {\em J. Combin. Theory Ser. B} {\bf 37} (1984), 96--100.
               
\bibitem{IKM} A.~Ili\'c, S.~Klav\v{z}ar, M.~Milanovi\'c, On distance-balanced graphs,
               {\em European J. Combin.} {\bf 31} (2010), 733--737. 

\bibitem{JKR} J.~Jerebic, S.~Klav\v zar, D.~F.~Rall, Distance-balanced graphs,  
                {\em Ann. Combin.} {\bf 12} (2008), 71--79.

\bibitem{KMMM} K.~Kutnar, A.~Malni\v c, D.~Maru\v si\v c, \v S.~Miklavi\v c,
     		Distance-balanced graphs: symmetry conditions, 
     		{\em Discrete Math.} {\bf 306} (2006),  1881--1894.

\bibitem{KMMM2} K.~Kutnar, A.~Malni\v c, D.~Maru\v si\v c, \v S.~Miklavi\v c,
     		The strongly distance-balanced property of the generalized Petersen graphs, 
     		{\em Ars Math. Contemp.} {\bf 2} (2009),  41--47.

\end{footnotesize}
\end{thebibliography}
\end{document}